\documentclass[a4paper,12pt]{article}

\usepackage{amsthm}
\usepackage{amsmath,amssymb,latexsym,amsfonts,mathrsfs}
\usepackage[dvipdfmx]{graphicx}

\usepackage{color}

\usepackage{fancyhdr}
\usepackage[top=30truemm, bottom=30truemm, left=25truemm, right=25truemm]{geometry}

\newcommand{\ep}{\varepsilon}
\newcommand{\nn}{\nonumber}

\newcommand{\SCR}[1]{{\mathscr #1}}

\newcommand{\D}[1]{{\mathscr D}( #1 )}

\theoremstyle{definition}
\newtheorem{Thm}{{\bf Theorem}}[section]

\newtheorem{Lem}[Thm]{{\bf Lemma}}

\newtheorem{Ass}[Thm]{{\bf Assumption}}

\newtheorem{Def}[Thm]{{\bf Definition}}
\newtheorem{Rem}[Thm]{{\bf Remark}}

\newcounter{Exami}

\newcommand{\Proof}[2][Proof]{
\begin{proof}[{\bf #1}]
#2
\end{proof}
}

\begin{document}

\begin{flushleft}
{\bf \Large Singularity for Solutions of Linearized KdV Equations
 } \\ \vspace{0.3cm} 
by {\bf \large Keiichi Kato $^{1}$ Masaki Kawamoto 
 $^{2}$ Koichiro Nanbu } \\ 
 $^{1}$ Department of mathematics, Faculty of science, Tokyo university
 of science, 1-3, Kagurazaka, Shinjuku-ku, Tokyo 162-8601, Japan. \\ 
 Email: kato@rs.tus.ac.jp \\ 
  $^{2}$ Department of Engineering for Production, Graduate School of Science and Engineering, Ehime University, 3 Bunkyo-cho Matsuyama, Ehime, 790-8577. Japan. \\ 
 Email: kawamoto.masaki.zs@ehime-u.ac.jp
\end{flushleft}

\begin{center}
\begin{minipage}[c]{400pt}
{\bf Abstract}. {\small 
We investigate the time propagation of singularity of a solution to linearized KdV equation by using the characterization of wave front sets via the wave packet transform (short time Fourier transform).   
}
\end{minipage}
\end{center}

\begin{flushleft}
\end{flushleft}
\section{Introduction}
In this paper, we investigate singularities of solutions to the linearized KdV equation written in the form 
\begin{align} \label{4}
\begin{cases}
& \partial _t u(t,x) + \partial _x^3 u(t,x) + \partial _x (a(t,x) u(t,x)) = 0, \\ 
& u(0,x) = u_0 (x) \in L^2({\bf R}_x), 
\end{cases}
\end{align}
where $\partial _t = \partial/\partial t$, $\partial _x = \partial /\partial x$, $(t,x)\in {\bf R}\times {\bf R}$, $a(t,x)$ a given smooth function which decays fast at infinity and $u(t,x)$ a real valued unknown function. \par
The motivation to consider the above problem is the following:
Consider the KdV equation
\begin{align} \label{1}
f_t(t,x) + a f_{}(t,x)f_x(t,x) + \gamma f_{xxx}(t,x) =0 , 
\end{align} 
where $(t,x)\in {\bf R}\times {\bf R}$ and $a$, $\gamma$ be positive parameters. 
For $x_0 \in {\bf R}$ and $abcd\gamma \neq 0$ with ratio $ac = 12 b^2 \gamma = 3d$, 
the following function which is called {\em soliton}
\begin{align} \label{2}
f(t,x) = c \cosh ^{-2} (b(x-dt-x_0))
\end{align}
solves the equation \eqref{1} with the initial condition $f(0,x) = c \cosh ^{-2} (b(x-x_0))$. Let $\ep >0$ be a sufficiently small constant and consider the nonlinear equation 
\begin{align} \label{3}
v_t + avv_x + \gamma v_{xxx} = \ep F(v).
\end{align}
Assuming that $v(t,x)$ can be written as the form $v(t,x) = f(t,x) + \ep w(t,x)$, then we have 
\begin{align*}
\ep \left( 
w_t + a(wf_x +f w_x) + \gamma w_{xxx}
\right) = \ep F(f) -\ep^2 w w_x - \ep^2 \int_0^1 F'(f+\ep \theta w)d\theta w. 
\end{align*}
If we neglect higher order term with respect to $\ep$ and smooth and rapidly decaying function $F(f)$, we have the first equation of \eqref{4} with 
$a=1$, $\gamma =1$, $b=1$ and $a(t,x) = f(t,x) $. 
Hence to investigate the behavior of a solution to the linearized problem \eqref{4} of the KdV equation would be a first step to know the singularities of 
the solutions to the KdV equation \eqref{3}. 
\par
Throughout this paper, we impose the following assumptions on the coefficient $a(t,x)$; 
 \begin{Ass} \label{A1}
Let $a \in C^1 ({\bf R} ; C^{\infty} ({\bf R})) $ be a real-valued function and suppose that there exists a real number $\rho >0$ such that 
for all $l _1 \in \{ 0,1\}$ and $ l_2\in {\bf N} \cup \{ 0\}$ we have 
 \begin{align} \label{17}
 \left| 
\partial _t^{l_1} \partial _x^{l_2} a(t,x) 
 \right|  \leq C_{l_1,l_2} (1+|x|)^{ -\rho - l_1-l_2}
 \end{align}
for some $C_{l_1,l_2} >0$. 
\end{Ass}
\begin{Rem}
Clearly the soliton $f$ in \eqref{2} satisfies Assumption \ref{A1}. 
\end{Rem}
Under this assumption, we have a unique solution of \eqref{4}. 
\begin{Thm} \label{T1}
Under the assumption \ref{A1}, there exists a unique solution $u= u(t,x)$ of \eqref{4}, which is included in $C({\bf R};L^2({\bf R}))$ and satisfies that for any $T \in {\bf R}$, 
\begin{align} \label{20}
\sup_{t \in [0,T]}\| u (t, \cdot ) \| _2 \leq C_T \| u_0\|_2. 
\end{align} 
Moreover if $u_0 \in H^3({\bf R})$, then $ u \in C^1({\bf R};L^2({\bf R})) \cap C({\bf R};H^3({\bf R}))$ and satisfies that for any $T \in {\bf R}$, 
\begin{align*}
\sup_{t \in [0,T]}\| u (t, \cdot)\| _{H^3 ({\bf R})} \leq C_T \| u_0\|_{H^3({\bf R})}. 
\end{align*} 
\end{Thm}
The aim of this paper is to consider the propagation of the singularities of the solution $u(t,x)$. More precisely we study the "wave front set" of $u(t,x)$, which indicates the place and the direction of the singularities of $u(t,\cdot)$. The notion of wave front set is introduced by L. H\"ormander in 1970 (\cite{Ho}).
\begin{Def}[Wave front sets]
For $f \in \SCR{S}' ({\bf R})$, the wave front set ${\rm WF}(f)$ of $f$ 
is a subset of ${\bf R}\times {\bf R}\backslash \{0\}$ which is determined as follows: 
We say $(x_0,\xi_0) \notin {\rm WF}(f)$ if there exists a function $\chi \in C_0^{\infty} ({\bf R})$ with $\chi (x_0) \neq 0$ and a conic neighborhood $\Gamma$ of $\xi _0$ such that for all $N \in {\bf N}$ there exists a positive constant $C_N$ so that 
\begin{align*}
\left| 
\widehat{\chi f} (\xi) \right|  \leq C_N (1+ |\xi|)^{-N}
\end{align*}   
holds for all $\xi \in \Gamma$, where $\hat{}$ denotes the Fourier transform. 
Otherwise we say $(x_0, \xi_0) \in {\rm WF}(f)$. 
\end{Def}

In order to analyze the wave front set of $u(t,\cdot)$, we employ the approach with using {\em Wave packet transform}, which was firstly considered by Folland \cite{Fo} and was developed by $\bar{\mathrm O}$kaji \cite{O}, Kato-Kobayashi-Ito\cite{IKK}, Pilipovic-Prangoski \cite{PP} and so on. Here the wave packet transform is defined as follows;
\begin{Def}[Wave packet transform]
Let $\varphi \in \SCR{S}({\bf R}) \backslash \{0\}$ and $f \in \SCR{S}'({\bf R})$. The wave packet transform $W_{\varphi} f(x, \xi)$ with window function $\varphi$ is defined as 
\begin{align} \label{5}
W_{\varphi} f(x, \xi) := \int_{{\bf R}} \overline{\varphi (y-x)} f(y) e^{-iy \xi} dy.
\end{align}
\end{Def} 

As is seen in \S{3}, the wave packet transform reduces \eqref{4} to the first order partial differential equations on ${\bf R} \times {\bf R}^{2}$ with error terms and such equations are closely related to {\em Hamilton equation}. By using the method of characteristics, we can obtain an integral equation which has the solutions to Hamilton equation associated to \eqref{4}. 
We can characterize the wave front sets of the solutions just by the asymptotic behavior of solution with respect to the parameter $\lambda$ which is equivalent to $|\xi|$ to the integral equation with the aide of characterization of wave front set via wave packet transform (Theorem 3.22 of \cite{Fo}, Theorem 2.2 of \cite{02} and Theorem 1.1 of \cite{IKK}). 
This is the merit of using the wave packet transform in the study of characterization of wave front set. We obtain the following characterization of wave front set;
\begin{Thm}\label{T2}
Let $u$ be a solution to \eqref{4} in $C({\bf R};L^2({\bf R}))$. Then the following two statements are equivalent: \\ ~~ \\ 
(i): $(x_0,\xi_0) \notin {\rm WF} (u(t_0, \cdot))$. \\ ~~ \\ 
(ii):  For any $N \in {\bf N}$, $b \geq 1$ and $\varphi_0 \in \SCR{S}({\bf R}) \backslash \{0\}$, there exist a neighborhood $K$ of $x_0$,  a neighborhood $\Gamma$ of $\xi_0$ and a constant $C_{N,b,\varphi_0} >0$ such that for all $x \in K$, $b^{-1} \leq |\xi| \leq b$ for $\xi \in \Gamma$ and $\lambda \geq 1$,  
\begin{align*}
\left| 
W_{\varphi_{\lambda} (-t_0) } u_0 (x(0;\lambda), \lambda \xi) 
\right|  \leq C_{N,b,\varphi} \lambda ^{-N}
\end{align*} 
holds, where $x(t;\lambda)$ is a solution to 
\begin{align*}
\dot{x}(t) &= -3 \lambda ^2 \xi ^2 + a(t,x(t)), \quad x(t_0) =x,  \\
\varphi_{\lambda} (t,x,\xi) &:= e^{-t(\partial _x^3 -3i \xi \partial _x^2) } \varphi_{0,\lambda} (x)
\end{align*} 
and $\varphi_{0, \lambda} (x) = \lambda ^{d/2} \varphi_0 (\lambda ^{d} x)$
 with $\min (\rho, 1/4) <d< 2\min (\rho, 1/4)$, i.e., 
\begin{align} \nn 
W_{\varphi_{\lambda} (-t_0) } u_0 (x(0;\lambda), \lambda \xi) 
&= 
\int_{\bf R} \overline{ \varphi_{\lambda} (-t_0,y-x(0;\lambda),\lambda \xi)} u_0(y) e^{-iy \cdot \lambda \xi} dy
\\ & = \int_{\bf R} 
\overline{e^{t_0(\partial_x^3-3i \lambda\xi \partial _x^2)}\varphi_{0,\lambda}
 ( y-x(0;\lambda))}  u_0 (y) e^{-iy \cdot \lambda \xi} dy. \label{15} 
\end{align}
\end{Thm}
The wave front sets of the solutions to Schr\"odinger equations have been studied in different way by Hassel-Wunsch \cite{H-W-1} and Nakamura \cite{Na-1}--\cite{Na-2}. 
\par
In the case that $a(t,x) \equiv 0$, the operator $x -3t \partial_x^3$ commutes 
with $\partial_t + \partial_x^3$, hence $(x -3t \partial_x^3)^l u(t,x)$ solves the equation if $u(t,x)$ solves the equation. This shows that 
$u(t,\cdot )$ is smooth for $t\ne 0$ if the initial data $u(0,x)= u_0(x)$ decays rapidly. 
Theorem \ref{T2} is a refinement of this phenomena. 
Our theorem determines the condition of initial data in which 
each point $(x_0, \xi_0)$ in ${\mathrm T}^*{\bf R}^n\backslash 0$ is not in the wave front set of the solution for given time $t_0$.  
$x(0,\lambda )= x + 3 \lambda^2 \xi^2 t_0$ in \eqref{15} corresponds to 
the operator $x -3t \partial_x^3$. If $\lambda $ tends to $\infty$, 
$x(0,\lambda )= x + 3 \lambda^2 \xi^2 t_0$ tends to $+\infty$, so 
the condition (ii) is the condition of the initial data $u_0(x)$ at $+\infty$. 
This phenomena is called {\em smoothing effect} and which has very important role in analyzing the solution to nonlinear equations. 
This kind of characterization has been studied in several equations. Beals \cite{B}, Rauch \cite{R} and Rauch-Reed \cite{RR1}-\cite{RR2} have studied for nonlinear hyperbolic equations.  Kato \cite{Ka2}, Doi \cite{Doi-1}--\cite{Doi-2}, Craig-Kappeler-Strauss \cite{C}, Biswas et al \cite{E}, Farah \cite{F},  Levandosky et al \cite{LSV} and Pokhozhaev \cite{P} have studied for dispersive equations including the KdV equation. 
\par
There are lots of works on wave packet transform, characterization of wave front set via wave packet transform and its applications. Cordoba-Fefferman \cite{C} have firstly introduced wave packet transform and applied it to hyperbolic equations. Folland \cite{Fo} has firstly characterized wave front set via wave packet transform.  $\bar{\mathrm O}$kaji \cite{O} and Kato-Kobayashi-Ito \cite{IKK} have improved the characterization of wave front set and have applied it to Schr\"odinger equations. Johansson \cite{J} and Pilipovic-Prangoski \cite{PP} have studied relationship between wave front set and short time Fourier transform, which is another name of wave packet transform. 
\par
As for the KdV equations, there are lots of studies on singularities of the solutions. 
In \cite{dBHK}, de Bouard, Hayashi and one of the authors have studied Gevrey regularizing effect for the KdV equation. In \cite{KO}, one of the authors and Ogawa have studied analyticity of solutions to the KdV equation with the initial data as Dirac's delta. Recent developments in studies for perturbed KdV equations can be seen in Biswas-Konar \cite{BK} and references their in. To the best of our knowledge, the wave front sets have not been considered for perturbed KdV equations.  Hence, as the first step for considering such issues, we consider the simplified linear model due to e.g., Mann \cite{Ma}. This attempt could be new and could be important not only mathematically but also physically. 
\section{Preliminaries and proof of Theorem \ref{T1}}
In this section, we introduce some notations and give the proof of Theorem \ref{T1}. Throughout this paper, $\left\| \cdot \right\|$ denotes the norm on $L^2({\bf R})$, and suppose that $x$ and $\xi$ are always included in the compact neighborhood of $x_0$ and conic neighborhood of $\xi _0$, respectively. For fixed $x_0$ and $t_0$, assume that $\lambda$ is always large enough compared to $|x|$ and $|t_0|$. $C$ will always denote a positive constant that does not depends on any parameters in consideration.  

\subsection{Proof of Theorem \ref{T1} }
Now we prove the existence and uniqueness of solution to \eqref{4}. In the proof, we employ the result of Enss-Veseli\'{c} \cite{EV}, that is
\begin{Lem}
 Let $H_0 (t)$ be a family of selfadjoint operators on $L^2({\bf R})$. Suppose that, for fixed $s \in {\bf R}$ and for every $t \in {\bf R}$, the operator
\begin{align*}
\frac{H_0(t)-H_0(s) }{t-s} \left( H_0(s) + i \right)^{-1} 
\end{align*} 
can be extended to bounded operators, then there uniquely exists a family of unitary operators $U_0(t,s)$ such that 
\begin{align*}
& i \frac{\partial }{ \partial t} U_0(t,s) = H_0(t) U_0(t,s), \quad i \frac{\partial }{ \partial s} U_0(t,s) = -U_0(t,s) H_0(s), \\ 
& U_0(t,s) = U_0 (t,0) U_0(s,0)^{\ast}, \quad U_0(s,s) = \mathrm{Id}_{L^2({\bf R}^n)}
\end{align*}
and 
\begin{align*}
U_0(t,s) \D{H_0(s)} \subset \D{H_0(s)}
\end{align*} 
hold. 
\end{Lem} 
Here we say $U_0(t,s)$ is the propagator for $H_0(t)$. Now we apply this lemma to our model. First, we reduce \eqref{4} to a form 
\begin{align*}
i\partial _t u(t,x) = -i \partial _x^3 u(t,x) + i\left(  (\partial _x (a(t,x)u(t,x) ) + a(t,x) (\partial _x u(t,x))   \right) /2 + i (\partial _x a) (t,x) u(t,x). 
\end{align*}
Let $s =0$. Employing the substitution $p = -i \nabla$, it follows 
\begin{align*}
i \partial _t u(t,x) = ( -p^3 - ( p a(t) + a(t) p)/2 + i(\partial _x a)(t) ) u(t,x), 
\end{align*}
where $a(t)$ is the multiplication operator of $a(t,x)$. Hence we define
\begin{align*}
H_0(t) = -p^3  - ( p a(t) + a(t) p)/2, \quad H(t) = H_0(t) + V(t) := H_0(t) + i (\partial _x a) (t). 
\end{align*} 
and first we prove the unique existence of unitary propagator for $H_0(t)$ and after we prove that $u(t,x)$ uniquely exists by using Duhamel's formula.

 It can be easily seen that $\mathscr{D}({H_0(0)}) = \mathscr{D}({p^3})$ holds since $a(t,x)$ is bounded, i.e., for all $\phi \in L^2({\bf R})$, 
 \begin{align*} 
 \sum_{k=0}^3 \left\| p^k (H_0(0) + i)^{-1} \phi \right\| \leq C \|  \phi \|
 .
 \end{align*} Then for $\phi \in L^2({\bf R})$
\begin{align*}
& \left\| 
\frac{(H_0(t)-H_0(0))}{t-0} (H_0(0)+i)^{-1}\phi
\right\|  \\ & = \left\| \left( 
p  \frac{a(t,x)-a(0,x)}{2t} + \frac{a(t,x)-a(0,x)}{2t} p
\right) (H_0(0) +i)^{-1} \phi
\right\| \\ & \leq 
\left\| 
\left( 
\frac{(\partial_x a) (t,x) - (\partial _x a) (0,x)}{2t}
\right) (H_0(0) + i)^{-1} \phi
\right\| + \left\| 
\frac{a(t,x) -a(0,x)}{t} p (H_0(0) +i)^{-1} \phi
\right\|  
\\ & \leq 
\left( 
C_{1,0} + \frac{C_{1,1}}{2}
\right) \left\| ( p + 1 )(H_0(0) + i)^{-1} \phi \right\|  \leq C \| \phi \| .
\end{align*}
holds. By this inequality and the result of \cite{EV}, we get the unique existence of $U_0(t,0)$. Now we show the unique existence of $U(t,0)$, the propagator for $H(t)$. By Duhamel's formula we have 
\begin{align*}
u(t,x) = U_0(t,0)u_0 + \int_0^t U_0(t,s) (\partial _x a (s)) u(s,x) ds. 
\end{align*}
Since $(\partial _x a (t))$ is bounded, we can use the standard argument in proving the existence of propagator, see e.g., Kato \cite{Ka}, IX \S{2}, and that proves the unique existence of the propagator $U(t,0) $ so that $u(t,x) =U(t,0)u_0$. Moreover since $\D{H_0(0)} = \D{H_0(s)} = \D{p^3} $ and $U_0(t,s) \D{H_0(s)} \subset \D{H_0(s) }$ hold, for $u(t) = u(t, \cdot)$, we have 
\begin{align*}
\left\| 
(H_0 (0 ) +i ) u(t)
\right\| &\leq \left\| (H_0(0) +i )U_0(t,0) u_0 \right\| + \int _0^t \Big\{
\left\| (H_0(0) +i) U_0(t,s) (H_0(s) + i)^{-1} \right\|_{\SCR{B}} \\ & \qquad \times\left\| (H_0(s) + i) (\partial _x a (s) ) (H_0(0) + i)^{-1} \right\| _{\SCR{B}} \left\| (H_0(0) + i) u(s) \right\| \Big\} ds, 
\end{align*}
where $\| \cdot \|_{\SCR{B}}$ denotes the operator norm from $L^2({\bf R})$ to itself. By the norm equivalence $c \| u \|_{H^3({\bf R})} \leq  \| (H_0(0) + i) u \| \leq C\| u \|_{H^3({\bf R})} $, we have, 
\begin{align*}
\left\| 
u(t)
\right\|_{H^3({\bf R})} \leq C \|  u_0 \|_{H^3({\bf R})} + C \int _0^t \left\| u(s) \right\|_{H^3({\bf R})} ds,  
\end{align*}
where the bound $\left\| (H_0(s) + i) (\partial _x a (s) ) (H_0(0) + i)^{-1} \right\| _{\SCR{B}} \leq C$ follows from $a (t , \cdot) \in C^{\infty} ({\bf R})$ and \eqref{17}.
By using this inequality, the proof of Theorem \ref{T1} completes. 

\section{Transformed equation via wave packet transform}
In this section, we transform \eqref{4} through wave packet transform and construct the solution of the transformed equation by using the solution of Hamilton-Jacobi equation.

Let us define $\varphi (t,x,\xi ) = e^{-t(\partial _x^3 + 3i \xi \partial _x^2)} \varphi_0 (x)$, $\varphi_0 \in \SCR{S}({\bf R}) \backslash \{ 0\}$. By \eqref{4}, 
\begin{align} \label{7}
W_{\varphi (t)} \partial _t u(t,x,\xi) +W_{\varphi (t)} \partial _x^3 u(t,x,\xi) +W_{\varphi (t)} \partial _x \left( 
au
\right)(t,x,\xi)=0
\end{align}
holds. Clearly 
\begin{align} \label{8}
W_{\varphi (t)} \partial _t u(t,x,\xi)  = \partial _t W_{\varphi (t)} u(t,x,\xi) -W_{\partial _t \varphi(t)} u(t,x,\xi)
\end{align}
holds by the definition of wave packet transform \eqref{5}. Furthermore, integration by parts and straightforward calculation shows 
\begin{align} \nn
W_{\varphi (t)} \partial _x^3 u(t,x,\xi) &= \int_{{\bf R}} \overline{\varphi (t,y-x,\xi)} \partial ^3 _y u(t,y) e^{-iy\xi} dy \\ &=  \nn 
 \int_{\bf R} \Big\{ -\overline{\partial _y^3 \varphi (t,y-x,\xi)} - \overline{3i \xi \partial _y^2 \varphi (t,y-x,\xi)}  \\ & \nn \quad + 3 \xi ^2\overline{\partial _y \varphi (t,y-x,\xi)} -i \xi ^3 \overline{\varphi(t,y-x, \xi)} 
\Big\}u(t,y) e^{-iy \xi}  dy \\&= \nn 
-W_{\partial _x^3 \varphi (t)} u(t,x,\xi) - W_{3i \xi \partial _x^2 \varphi (t)}u(t,x,\xi) \\ & \qquad -3 \xi ^2 \partial _x W_{\varphi (t)} u(t,x,\xi) -i \xi ^3 W_{\varphi (t)} u(t,x,\xi). \label{6}
\end{align}
By the equation 
\begin{align*}
W_{\partial _t \varphi (t) } u(t,x,\xi) +  W_{\partial _x^3 \varphi (t) }u(t,x,\xi) + W_{3i \xi \partial _x^2 \varphi (t) } u(t,x,\xi) = 0 ,
\end{align*} 
and equations \eqref{7}, \eqref{8} and \eqref{6}, we get 
\begin{align*}
& \partial _t W_{\varphi (t)}u(t,x,\xi) \\ &= W_{\varphi(t)} \partial _t u(t,x,\xi) + W_{\partial _t \varphi (t)} u(t,x , \xi) \\ &= 
 3 \xi ^2 \partial _x W_{\varphi (t)}u(t,x, \xi) + i \xi ^3 W_{\varphi (t)} u(t,x, \xi) - W_{\varphi (t)} \partial _x(au)(t,x, \xi). 
\end{align*}
To deal with the term $ W_{\varphi (t)} \partial _x(au)(t,x, \xi)$, we use the Taylor expansion; 
\begin{align*}
a(t,y) &= a(t,x+(y-x)) = a(t,x) + \sum_{k=1}^{L-1} \frac{\partial _x^k a(t,x)}{k !} (y-x)^k + r_L(t,x,y) (y-x)^L, \\ &=: a(t,x) + \tilde{r} (t,x,y) \\
 r_L(t,x,y) & := \frac{1}{(L-1) !} \int_0^1 \partial _x^L a(t,x + \theta (y-x) ) (1- \theta)^{L-1} d \theta .
\end{align*}
By employing similar calculations as in \eqref{6}, 
\begin{align*}
& W_{\varphi(t)} (\partial _x(au))(t,x,\xi) \\ & =- \int_{\bf R} \overline{\partial _y \varphi (t,y-x,\xi)} (au) (t,y) e^{-iy \xi} dy + i \xi \int_{\bf R} \overline{\varphi (t,y-x,\xi)} (au) (t,y) e^{-iy \xi} dy \\ &= 
-\int_{\bf R} \overline{\partial _y \varphi (t,y-x,\xi)} \left( a(t,x) + \tilde{r} (t,x,y) \right)u (t,x, \xi) dy \\ & \qquad + i \xi \int_{\bf R} \overline{\varphi (t,y-x,\xi)} (au) (t,y) e^{-iy \xi} dy .
\end{align*}
Hence we finally obtain the equations 
\begin{align} \label{9}
\begin{cases}
\left( \partial _t + \left( 
-3 \xi ^2 + a(t,x) \right) \partial _x -i \xi ^3 +i \xi a(t,x)
\right) W_{\varphi (t)} u(t,x,\xi) =Ru(t,x,\xi), \\ 
W_{\varphi (0)} u(0,x,\xi) = W_{\varphi_0} u_0(x,\xi), 
\end{cases}
\end{align}
where 
\begin{align*}
Ru(t,x,\xi) = \sum_{k=1}^{L-1} \frac{\partial _x^k a(t,x)}{k!} \left( 
-W_{x^k \partial _x \varphi (t)} u(t,x,\xi) + i \xi W_{x^k \varphi (t)} u(t,x,\xi)
\right) + R_L u(t,x,\xi)
\end{align*}
with 
\begin{align*}
R_L u(t,x,\xi) &= - \int \overline{\partial _y \varphi (e,y-x,\xi)} r_L (t,x,y) (y-x)^L u(t,y) e^{-i\xi y} dy \\ & \quad + i \xi \int \overline{\varphi (t,y-x,\xi)} r_L (t,x,y) (y-x)^L u(t,y) e^{-i \xi y} dy. 
\end{align*}
Now we solve \eqref{9} by using the method of characteristic curve. Fix $t_0 >0$ and let $x(t)$ and $\xi (t) (\equiv \xi)$ be solutions to 
\begin{align} \label{10}
\begin{cases}
\dot{x} (t) &= -3 \xi ^2 + a(t,x(t)) , \\ 
\dot{\xi} (t) &= 0, 
\end{cases}
\qquad 
\begin{cases}
x(t_0) =x, & \\ 
\xi (t_0) = \xi. &
\end{cases} 
\end{align}
By employing similar technique as in \S{5} of \cite{IKK}, we get 
\begin{align*}
W_{\varphi (t)}u(t,x(t),\xi) &:= e^{i\int_0^t \left( \xi ^3 - \xi a(\tau,x(\tau)) d\tau \right) } W_{\varphi_0 } u_0(x(0), \xi) \\ & \qquad  + \int_0^t e^{i\int_s^t \left( \xi ^3 - \xi a(\tau , x (\tau)) \right) d \tau } 
Ru(s,x(s),\xi) ds  
\end{align*}
will be a solution to \eqref{9}. Since the solution $W_{\varphi (t)} u(t,x(t),\xi)$ does not depend on the choice of $\varphi $, we substitute $\varphi_{\lambda} (t-t_0,x,\xi)$ instead of $\varphi (t,x,\xi)$, where we remark that 
\begin{align*}
\varphi_{\lambda} (t-t_0,x,\xi) = e^{-(t-t_0)\left( \partial_x^3 - 3i \xi \partial_x^2  \right) } \varphi_{0, \lambda}, \qquad \varphi_{0,\lambda}(x) = \lambda ^{d/2} \varphi_0 (\lambda ^{d}x )
\end{align*} 
with $0<d\le 1/2$, which is fixed later. 
Then we have 
\begin{align*}
W_{\varphi_{\lambda} (t-t_0) } u(t,x(t),\xi) &= e^{i\int_0^t \left( \xi ^3 - \xi a(\tau,x(\tau) )  \right) d \tau } W_{\varphi_{\lambda}(-t_0) } u_0(x(0),\xi) \\ 
& \qquad + \int_0^t e^{i \int_s^{t} \left( 
\xi ^3 -\xi a(\tau,x(\tau))
\right) d \tau} Ru (s,x(s),\xi  ) ds. 
\end{align*}
Furthermore, by introducing $x(t;\lambda)$ as solution to 
\begin{align*}
\dot{x} (t) = -3 \lambda^2 \xi ^2 +a (t,x(t)), \qquad x(t_0) = x,  
\end{align*}
we finally get 
\begin{align} \nn 
W_{\varphi_{\lambda} (t-t_0) } u(t,x(t;\lambda),\lambda \xi) &= e^{i\int_0^t \left( \lambda ^3\xi ^3 - \lambda \xi a(\tau,x(\tau;\lambda) )  \right) d \tau } W_{\varphi_{\lambda}(-t_0) } u_0(x(0;\lambda),\lambda \xi)\label{} \\ 
& \qquad + \int_0^t e^{i \int_s^{t} \left( 
\lambda ^3\xi ^3 - \lambda \xi a(\tau,x(\tau;\lambda))
\right) d \tau} Ru (s,x(s;\lambda), \lambda\xi  ) ds. \label{13}
\end{align}
By using this solution, we prove Theorem \ref{T2}.

\section{Proof of Theorem \ref{T2}} 
In this section, we prove Theorem \ref{T2}. The approach is based on the argument of \cite{IK}. The proof completes by showing both relations $(i) \Rightarrow (ii)$ and $(ii) \Rightarrow (i)$. The second relation can be proven quite simple by using the approach of proving $(ii) \Rightarrow (i)$ and hence we first prove that $(ii) \Rightarrow (i)$ and after we prove $(i) \Rightarrow (ii)$. 

Before the proof of main theorem, we introduce the following important lemma; 
\begin{Lem} \label{L1}
Let $0<\theta <2 $. For some $\lambda_0 \ge 1$, 
the solution $x(s;\lambda)$ of \eqref{10} satisfies 
\begin{align*}
\left| 
x(s;\lambda)
\right|  \geq \frac{3}{2b^2} \lambda ^2 |s-t_0 | 
\end{align*}
 for all $x \in K$, $\xi \in \Gamma$ with $b^{-1} \leq |\xi| \leq b$, $\lambda \geq \lambda_0$ and $s\ge 0$ with $|t_0 -s| \geq \lambda^{-\theta} $. 
\end{Lem}
\Proof{
The proof is almost the same as the proof of (11) in \cite{IK}, hence we only give the sketch of the proof. We use Picard's iteration method. Let us define 
\begin{align*}
\begin{cases}
x^{(N+1)}(s) &:= x -3 (s-t_0) \lambda ^2 \xi^2 - \int_{t_0}^s (s-s_1) \nabla _x V(s_1, x^{(N)} (s_1)) d s_1 ,  \\
x^{(0)} (s) &:= x -3 (s-t_0) \lambda ^2 \xi^2 
\end{cases}
\end{align*}
and use the induction scheme; Since $\theta <2$, by taking $\lambda >1$ large enough, it can be seen that    
\begin{align*}
|x^{(0)} (s)| \geq 3|s-t_0| \lambda ^2 \xi ^2 - |x| \geq 3b^{-2}|s-t_0| \lambda ^2  - C \geq 3b^{-2}|s-t_0| \lambda ^2 /2 
\end{align*}
holds, where we use $K$ is compact. Hence suppose that 
for some $M \in {\bf N} \cup \{ 0 \} $, 
\begin{align*}
|x^{(M)} (s)| \geq 3b^{-2}|s-t_0| \lambda ^2 /2 
\end{align*}
holds. Then we prove that 
\begin{align*}
|x^{(M+1)} (s)|  \geq 3b^{-2}|s-t_0| \lambda ^2 /2 
\end{align*}
also holds. Indeed 
\begin{align*}
|x^{(M+1)} (s)| \geq 5 b^{-2} |s-t_0|  \lambda ^2 /2  - \left| 
\int_{t_0} ^s (s-s_1) \nabla _x V(s_1, x^{(N)} (s_1)) ds_1
\right|  
\end{align*}
holds and by using Assumption \ref{A1}, for a positive constant $c>0$, we have    
\begin{align*}
|x^{(M+1)} (s)| & \geq 5b^{-2}|s-t_0| \lambda ^2 /2 -C  \int_{t_0}^{s} |s-s_1| \left( 
1+ c \lambda ^2 b^{-2} |(s_1-t_0)|
\right) ^{-\rho -1}ds_1 \\ & \geq 
5b^{-2}|s-t_0| \lambda ^2 /2   -C  \\ & \geq 
3b^{-2} |s-t_0| \lambda ^2 /2, 
\end{align*}
which is the desired result.
}
\subsection{Proof for  $(ii) \Rightarrow (i)$}
It suffices to prove the following assertion $P(N, \varphi_0)$ for all $N \in {\bf N}$ and $\varphi _0 \in \SCR{S}({\bf R}) \backslash \{ 0 \}$ by induction. \\ ~~ \\
{\bf $P(N,\varphi_0)$:} For $b \geq 1$, there exist a neighborhood $K$ of $x_0$, a neighborhood $\Gamma$ of $\xi_0$ and $C_{N,b, \varphi_0}$ such that for all $x \in K$, $\xi \in \Gamma$ with $b^{-1} \leq |\xi| \leq b$, $\lambda \geq 1$ and $0 \leq t \leq t_0$,  
\begin{align*}
\left| 
W_{\varphi_{\lambda}(t-t_0) } u(t,x(t;\lambda), \lambda \xi)
\right| \leq C_{N,b,\varphi_0} \lambda^{-N}
\end{align*}
holds. \\ ~~ \\
In fact, if the assertion $P(N, \varphi_0)$ holds for $N \in {\bf N}$ and $\varphi _0 \in \SCR{S}({\bf R}) \backslash \{ 0 \}$,   
\begin{align} \label{16}
\left| 
W_{\varphi_{0, \lambda }} u(t_0,x, \lambda \xi)
\right| \leq C_{N,b,\varphi_0} \lambda^{-N}
\end{align}
holds for all $x \in K$, $\xi \in \Gamma$ with $b^{-1} \leq |\xi| \leq b$, $\lambda \geq 1$ by taking $t=t_0$, since $\varphi_{\lambda} (0) = \varphi_{0,\lambda}$ and $x(t_0;\lambda) =x$. Then Theorem 1.1 of \cite{IKK} which is a refinement of Theorem 2.2 of \cite{02} and Theorem 3.22 of \cite{Fo}, shows that \eqref{16} is equivalent to $(x_0,\xi_0) \notin {\rm WF}(u(t_0, \cdot))$, i.e., we have (i) in Theorem \ref{T2}. Hence we prove $P(M M_0, \varphi_0)$ by the induction with respect to $N$ for some small constant $M_0$. We note that we can assume that $0 < \rho \leq 1/4$ without loss of generality; \\ ~~ \par
(I): Proof of  $P(0, \varphi_0)$. \\ 
Since $u(t,x)$ is in $C({\bf R};L^2({\bf R}))$, 
we have by Schwarz's inequality, $L^2$ conservativity for $e^{-t(\partial_x^3 -3i \xi \partial _x^2)}$ and \eqref{20} 
\begin{equation}
 \left| 
W_{\varphi_{\lambda}(t-t_0) } u(t,x(t;\lambda), \lambda \xi)
\right| 
\le C \Vert \varphi_{\lambda}(t-t_0)\Vert_{L^2} \Vert u(t,\cdot)\Vert_{L^2}
\leq \Vert \varphi_0\Vert_{L^2} \Vert u_0\Vert_{L^2}. 
\end{equation}
\par
(II): Assuming that for $M\in {\bf N}$, $P(MM_0, \varphi _0)$ is valid for all $\varphi _0 \in \SCR{S}({\bf R}) \backslash \{ 0 \}$, we show that $P((M+1)M_0, \varphi_0)$ holds for all $\varphi _0 \in \SCR{S}({\bf R}) \backslash \{ 0 \}$ for some small constant $0< M_0 < 1/2$. Since the condition (ii) in Theorem \ref{T2} says that the modulus of the first term of the right hand side of \eqref{13} is less than or equal to $C \lambda ^{-(M+1)M_0}$, it is enough to prove that
\begin{align*}
\left| 
\int_0^t Ru(s,x(s;\lambda),\lambda \xi) ds 
\right| \leq C_{N,b,\varphi_0} \lambda ^{-(M+1)M_0}
\end{align*}
holds. In the followings, for simplicity, we may use notations $x_{\lambda}(t) = x(t;\lambda)$, $U(t,\xi) = e^{-t(\partial_x^3 -3i \xi \partial _x^2)}$ and $(Uf) (t,x) = U(t,\xi) f(x)$. We put 
\begin{equation*}
Ru(t,x,\xi) = \Gamma_1(t) + R_L u(t,x,\xi), 
\end{equation*}
where
\begin{equation*}
\Gamma_1(t) = \sum_{k=1}^{L-1} \frac{\partial _x^k a(t,x)}{k!} \left( 
-W_{x^k \partial _x \varphi (t)} u(t,x,\xi) + i \xi W_{x^k \varphi (t)} u(t,x,\xi)
\right) . 
\end{equation*}
\par
(II-I) Estimation for $ \Gamma _1(t) $.
By the simple calculation, we can see 
\begin{align*}
(x^k) \varphi_{\lambda} (t,x,\lambda \xi) & = U(t, \lambda \xi) (x+3t \partial _x^2-6it \lambda \xi \partial _x  ) ^k \varphi_{0, \lambda} (x) \\ &= 
U(t, \lambda \xi) \sum_{\alpha + \beta + \gamma = k} C_{\alpha , \beta, \gamma} \lambda ^{\gamma} \xi ^{\gamma} t^{\beta + \gamma} (x^{\alpha}\partial _x^{2 \beta + \gamma} ) \varphi_{0, \lambda} (x) \\ &= 
\sum_{\alpha + \beta + \gamma = k} C_{\alpha , \beta, \gamma} \lambda ^{-d\alpha + 2d\beta +(1+d)\gamma } \xi ^{\gamma} t^{\beta + \gamma} U(t,\lambda ,\xi)(x^{\alpha}\partial _x^{2 \beta + \gamma}  \varphi_{0}) _{\lambda} (x) . 
\end{align*}
By the inductive hypothesis, we have
\begin{align*}
&\left| \int_0^t \Gamma _1 (s) ds \right| \\  & \leq C 
\int_0^t \sum_{k=1}^{L-1} \sum_{\alpha + \beta + \gamma = k} \left| 
\partial _x^k a(s,x_{\lambda}(s)  )
\right| \lambda ^{-d\alpha  + 2d\beta + (1+d)\gamma } |\xi|^{\gamma} |s-t_0|^{\beta + \gamma}  \\ & \quad \times 
\left( 
\lambda ^{d} \left| 
W_{U(x^k\partial _x \varphi_0)_{\lambda}(s-t_0) } u(s,x_{\lambda}(s), \lambda \xi )
\right|  + \lambda |\xi| \left| W_{U(x^k \varphi_0)_{\lambda} (s-t_0)} u(s,x_{\lambda}(s), \lambda \xi ) \right| 
\right) ds \\ & \leq 
C \int_0^t \sum_{k=1}^{L-1} \sum_{\alpha + \beta + \gamma =k} \left( 
1+ |x_{\lambda} (s)| \right) ^{-\rho -k} \lambda ^{-d\alpha + 2d\beta +(1+d) \gamma } |s-t_0|^{\beta + \gamma} \lambda \lambda ^{-MM_0} ds \\ 
& \leq \Gamma _2 + \Gamma _3
\end{align*}
with, for $0 < \theta <2$, 
\begin{align*}
\Gamma _2 := C \int_{|s-t_0| \geq \lambda ^{- \theta}} \sum_{k=1}^{L-1} \sum_{\alpha + \beta + \gamma =k} \left( 
1+ |x_{\lambda} (s)| \right) ^{-\rho -k} \lambda ^{-d\alpha + 2d\beta +(1+d) \gamma +1 } |s-t_0|^{\beta + \gamma}  \lambda ^{-MM_0} ds 
\end{align*}
and 
\begin{align*}
\Gamma _3 := C \int_{|s-t_0| \leq \lambda ^{- \theta}} \sum_{k=1}^{L-1} \sum_{\alpha + \beta + \gamma =k} \left( 
1+ |x_{\lambda} (s)| \right) ^{-\rho -k} \lambda ^{-d\alpha + 2d\beta +(1+d) \gamma +1 } |s-t_0|^{\beta + \gamma}  \lambda ^{-MM_0} ds .
\end{align*}
By using Lemma \ref{L1}, we have 
\begin{align*}
\Gamma_2 & \le C \int_{|s-t_0| \geq \lambda ^{- \theta}} \Big\{  \sum_{k=1}^{L-1} \sum_{\alpha + \beta + \gamma =k}  \left( 
1+ 3b^{-2}\lambda^2 |s-t_0| /2 \right) ^{-\rho -k} \\ &  \phantom{xxxxxxxxxxxxxxx} \times \lambda ^{-d\alpha + 2d\beta +(1+d) \gamma +1 }  |s-t_0|^{\beta + \gamma} \lambda ^{-MM_0} \Big\} ds .
\end{align*}
Here by employing the inequality 
\begin{align*}
k \geq \mathrm{min} (k, \beta + \gamma +1) -2d \min (k- \beta - \gamma ,1), 
\end{align*}  
we have the estimate 
\begin{align*}
\Gamma _2  & \leq   \sum_{k=1}^{L-1} \sum_{ \alpha + \beta + \gamma =k  }
\lambda ^{-d\alpha + 2d\beta +(1+d) \gamma +1 -2\min (k,\beta +\gamma +1 )-2\rho   +2d\min(k-\beta-\gamma, 1) -MM_0}  \\
&\phantom{xxxxxxxxxxxxxxx}\times
\int_0^t |s-t_0|^{\beta + \gamma -\min (k,\beta +\gamma +1)-\rho +d\min(k-\beta-\gamma, 1) }  ds \\
&\leq C
\lambda ^{d - 2\rho -MM_0},
\end{align*}
since $-d\alpha + 2d\beta +(1+d) \gamma  -2\min (k,\beta +\gamma +1)+2d\min(k-\beta-\gamma, 1)$ takes the maximum $d-1$ when $k=1$, $\alpha =\beta =0$ and $\gamma =1$. 
If we take $\rho <d<2\rho $ and $0< M_0 \le 2\rho -d$, we have 
$$
\Gamma _2 \le C\lambda^{-(M+1)M_0}. 
$$
Since $\int_{|s-t_0| \leq \lambda ^{- \theta}} ds \leq C \lambda ^{-\theta}$ holds, by taking $\theta$ so that $1 + d < \theta <2$, we get 
\begin{align*}
\Gamma _3 \leq C \lambda^{-\theta - d \alpha + 2d \beta + (1+d) \gamma +1- \theta (\beta + \gamma)  -MM_0} \leq C \lambda^{-(M+1) M_0}
\end{align*} 
for $0< M_0 <d$. Consequently we get 
\begin{align*}
\left| \int_0^t \Gamma _1 (s) ds \right| \leq C \lambda ^{-(M+1) M_0}.
\end{align*}
\par
(II-II) Estimation for $R_L u(t,x(t;\lambda), \lambda \xi)$; \\ 
Let us define $\psi_1 \in C_0^{\infty} ({\bf R}) $ and $\psi _2 \in C^{\infty} ({\bf R})$ as follows
\begin{align*}
\psi_1 (x) = 
\begin{cases}
1  & |x| \leq 1, \\ 
0 & |x| \geq 2, 
\end{cases} 
\qquad 
\psi_2 (x) = \begin{cases}
0  & |x| \leq 1, \\ 
1  & |x| \geq 2, 
\end{cases} 
\qquad \psi_1 + \psi_2 =1. 
\end{align*}
Devide $R_L u = R_{L,1}u + R_{L,2} u $ with 
\begin{align*}
& R_{L,j} u(s,x(s;\lambda), \lambda \xi) \\ & 
= \int 
\overline{
\left( 
\partial _x(\varphi_{\lambda}) +i \lambda \xi \varphi_{\lambda}
\right) (s, y-x(s;\lambda) ,\lambda \xi) 
} r_L (s, x(s; \lambda) ,y ) (y-x(s;\lambda))^L \\ & \qquad \times 
\psi_j \left( 
\frac{\lambda^{1/4}(y-x(s;\lambda)) }{1+ \lambda ^2 |s-t_0|}
\right) u(y) e^{-iy \lambda \xi} dy.
\end{align*} 
First, we estimate the term associated to $R_{L,1}$. \\ 
$\cdot$ The case that  $|s-t_0| \leq \lambda^{- 15/8}$. \\
On the support of $\psi _1$ and $|s-t_0| \leq \lambda^{- 15/8}$, it holds that 
\begin{align*}
|y-x(s;\lambda)| \leq 2 \lambda^{-1/4} (1+ \lambda^2 |s-t_0| ) \leq C \lambda^{-1/8}, 
\end{align*}
and which deduces 
\begin{align*}
& \int_{|s-t_0| \leq \lambda^{- 15/8}} \left| 
R_{L,1} u(s,x_{\lambda } (s),\lambda \xi)
\right| ds  \\ & \leq C \int_{|s-t_0| \leq \lambda^{- 15/8}} \int \left( 
\lambda^{d} \left| 
(\partial _x \varphi)_{\lambda} (s,y-x_{\lambda} (s), \lambda \xi) 
\right| + \lambda |\xi| \left| 
\varphi_{\lambda} (s,y-x_{\lambda}(s), \lambda \xi )
\right| 
\right) \\ & \qquad \times 
\left| 
r_L (s,x_{\lambda}(s),y ) 
\right| \left| y-x \right|^L |u(y)| dy ds \\ & \leq 
C \lambda ^{-L/8 +1 -15/8 }  \left( 
\left\| 
(\partial _x \varphi) _{\lambda}
\right\| + \left\| 
\varphi_{\lambda}
\right\| 
\right)  \left\| u \right\|. 
\end{align*} 
$\cdot$ The case of  $|s-t_0| \geq \lambda^{-15/8}$; \\ 
On the $\mathrm{supp} \left\{ \psi _1 \right\} (\lambda ^{1/4}(y-x(s;\lambda)) /(1+ \lambda ^2 |s-t_0|) )$, it holds that 
\begin{align*}
\left| 
y-x(s; \lambda)
\right| \leq 2 \lambda^{-1/4} (1+ \lambda ^2|s-t_0| ).
\end{align*} 
Together with Lemma \ref{L1}, we get for $\lambda \geq \lambda _0$
\begin{align*}
\left| 
x(s;\lambda) + \theta (y-x(s;\lambda))
\right|  \geq |x(s; \lambda)| - |y-x(s; \lambda)| \geq C \lambda ^2 |s-t_0|.
\end{align*}
Consequently, 
\begin{align*}
\left| 
\partial _x^L a(s,x(s;\lambda) + \theta (y-x(s;\lambda)))
\right| &\leq C \left( 
1+ |x(s;\lambda) + \theta (y-x(s;\lambda)) |
\right) ^{-\rho -L} \\ & \leq C(1+ \lambda ^2|s-t_0| )^{-\rho-L}
\end{align*}
holds. Then we have 
\begin{align*}
& \int_{|s-t_0| \geq \lambda ^{- \theta}} \left| 
R_{L,1} u(s,x_{\lambda } (s),\lambda \xi)
\right| ds  \\ & \leq C \int_{|s-t_0| \geq \lambda ^{- \theta}} \int 
\left( 
\lambda^{d} \left| 
(\partial_x \varphi)_{\lambda} (s,y-x_{\lambda} (s), \lambda \xi )
\right| + \lambda |\xi| \left| \varphi_{\lambda} (s,y-x_{\lambda} (s), \lambda \xi) \right|
\right) \\ & \qquad \times 
\left( \left(
1+ \lambda ^2 |s-t_0| \right)^{-\rho -L} \lambda ^{-L/4} \left( 
1+ \lambda ^2 |s-t_0|
\right) ^L |u(y)| dy ds  
\right) \\ & \leq C \lambda ^{-L/4 +1} \int_0^t \left( 
1+ \lambda^2 |s-t_0| 
\right)^{- \rho}  ds \times \left( 
\left\| 
(\partial _x \varphi) _{\lambda}
\right\| + \left\| 
\varphi_{\lambda}
\right\|
\right)  \left\| u \right\| \\ & \leq 
C \lambda ^{-L/4-1}.
\end{align*}
By choosing $L$ large, we have 
\begin{align*}
\int_0^t \left| 
R_{L,1} u(s,x_{\lambda } (s),\lambda \xi)
\right| ds  \leq 
C \lambda ^{-(M+1)M_0}.
\end{align*}
Next, we estimate the term associated to $R_{L,2}$. For $m \in {\bf N}$, 
\begin{align*}
\left( 
1+|x|^2
\right)^{m} \varphi_{\lambda} (t,x,\lambda \xi) &= e^{-t(\partial _x^3 -3i \lambda \xi \partial _x^2)} \left\{ 
1+ \left( 
x+3t\partial _x ^2 -6it \lambda \xi \partial _x
\right) 
\right\} ^m \varphi_{0,\lambda} (x) \\ &= 
e^{-t(\partial _x^3 -3i \lambda \xi \partial _x^2)} \sum_{\alpha + \beta + \gamma \leq 2m} C_{\alpha , \beta, \gamma} \lambda ^{\gamma} \xi ^{\gamma} t^{\beta + \gamma} \left( 
x^{\alpha} \partial _x^{2 \beta + \gamma}
\right) \varphi_{0, \lambda} (x) \\ &=
e^{-t(\partial _x^3 -3i \lambda \xi \partial _x^2)} \sum_{\alpha + \beta + \gamma \leq 2m} C_{\alpha , \beta, \gamma} \lambda ^{- d\alpha  + 2d\beta + (1+d) \gamma} \xi^{\gamma} t^{\beta + \gamma} \left( x^{\alpha} \partial _x^{2 \beta + \gamma} \varphi_0 \right) (x) . 
\end{align*}
By using this, we can obtain 
\begin{align*}
& \int_0^t \left| 
R_{L,2} u(s,x_{\lambda} (s), \lambda \xi )
\right| ds  \\ & 
\leq \int_0^t \int 
\left( 
1+ |y-x_{\lambda}(s) | ^2
\right)^{-m}\left( 
1+ |y-x_{\lambda}(s) | ^2
\right)^{m} \left( \left| 
\partial _x (\varphi_{\lambda})  \right|  + \left| i \lambda \xi \varphi_{\lambda}  
(s-t_0, y-x_{\lambda }(s)  , \lambda \xi )
\right|  \right) \\ & \qquad \times 
\left| r_L (s,x_{\lambda} (s),y) \right| \left| y-x_{\lambda}(s) 
\right|^L |u(y)| dy ds \\ & 
\leq C 
\int_0^t \int \left(1+ |y-x_{\lambda}(s) |^2 \right)^{-m} \\ & \quad \times 
\Big( 
\lambda ^{d} \sum_{\alpha + \beta + \gamma = 2m} \lambda^{- d\alpha  + 2d \beta +(1+d) \gamma } |s-t_0|^{\beta + \gamma} \left| U(s-t_0) \left( x^{\alpha} \partial_x^{2 \beta + \gamma +1}      
 \varphi _0 \right)_{\lambda}(y-x) \right| \\ & \qquad + 
 \lambda \sum_{\alpha + \beta + \gamma = 2m} \lambda ^{-d\alpha  + 2d\beta +(1+d) \gamma } |s-t_0|^{\beta + \gamma} \left| 
  U(s-t_0) \left( x^{\alpha} \partial_x^{2 \beta + \gamma }      
 \varphi _0 \right)_{\lambda}(y-x)
 \right| 
\Big)  \\ & \qquad \quad \times \left|y-x_{\lambda} (s) \right|^L |u(y)| dyds \\ & \leq 
C \int_0^t \int \left( 
1+ |y-x_{\lambda}(s) |^2 
\right) ^{-m + L/2} \lambda^{2(1+d)m+1} (1+|s-t_0|^{2m}) \\
& \phantom{xxxxxxxxxxxxxx}
\times\left| 
U(s-t_0) \left( \partial _x^{2m} \varphi_0
\right)_{\lambda} (y-x) 
\right| \left| u(y) \right| dy ds .
\end{align*}
Since $|y-x(s;\lambda)| \ge \lambda^{-1/4} (1+ \lambda^2 |s-t_0| )$ holds in the support of $\psi _2$, 
Schwarz's inequality shows that 
\begin{align*}
& \int_0^t \left| 
R_{L,2} u(s,x_{\lambda} (s), \lambda \xi )
\right| ds  \\
 & \leq 
C \int_0^t \int \left( 
\lambda^{-1/4} \lambda^2 |s-t_0| 
\right) ^{-2m + L} \lambda^{2(1+d)m+1} (1+|s-t_0|^{2m}) \\
&\phantom{xxxxxxxxxxx}
\times\left| 
U(s-t_0) \left( \partial _x^{2m} \varphi_0
\right)_{\lambda} (y-x) 
\right| \left| u(s, y) \right| dy ds \\
 & \leq 
C \lambda^{-(3/2-2d)m +1+(7/4)L}
\int_0^t 
(1+|s-t_0|^{L} )
\left\| U(s-t_0) \left( \partial _x^{2m} \varphi_0 \right)_{\lambda} (y-x) \right\| \left\| u(s, y)\right\|
ds \\
 & \leq 
C \lambda^{-(3/2-2d)m +1+(7/4)L}
\int_0^t 
(1+|s-t_0|^{L}) ds
\left\| \left( \partial _x^{2m} \varphi_0 \right)_{\lambda} (y-x) \right\| 
\left\| u_0(\cdot )\right\|\\
& \leq C \lambda ^{-(1/2)m +1+(7/4)L}
\leq C \lambda ^{-(M+1)M_0}, 
\end{align*}
if we take $m$ sufficiently large compared to $L$. 
Hence we get $P((M+1)M_0, \varphi_0)$  for any $\varphi_0 \in \SCR{S}({\bf R}) \backslash \{ 0 \}$. This completes the proof. \qed

\subsection{Proof for $(i) \Rightarrow (ii)$}
If we transform the equation \eqref{9} with the initial condition 
$W_{\varphi(-t_0)}u(t_0,x,\xi )$ for $t=t_0$, we have 
\begin{align*}
W_{\varphi_{\lambda} (t-t_0) } u(t,x_{\lambda }(t),\lambda \xi) &= e^{i\int^t_{t_0} \left( \xi ^3 - \xi a(\tau,x(\tau) )  \right) d \tau } W_{\varphi_{0,\lambda} } u(t_0,x_{\lambda}(t_0),\lambda \xi) \\ 
& \qquad + \int_{t_0}^t e^{i \int_{t_0}^{s} \left( 
\xi ^3 -\xi a(\tau,x_{\lambda}(\tau))
\right) d \tau} Ru (s,x_{\lambda }(s),\lambda \xi  ) ds, 
\end{align*}
where $x_{\lambda}(t)$ is a solution of $\dot{x}(t) = -3 \lambda ^2 \xi ^2 + a(t,x(t))$ with $x(t_0) =x$. 
Hence we have 
\begin{equation*}
\left| 
W_{\varphi _{\lambda} (t-t_0) } u_0 (x_{\lambda} (t), \lambda \xi)
\right| \leq \left| 
W_{\varphi _{0,\lambda}} u(t_0,x_{\lambda}(t), \lambda \xi )
\right| + \int_t^{t_0} \left| 
Ru(s,x_{\lambda}(s), \lambda \xi )
\right|  ds. 
\end{equation*}
We note that the condition $(i)$ is equivalent to $|W_{\varphi_{0,\lambda}} u(t_0,x, \lambda \xi)| \leq C \lambda ^{-N}$. Hence by the same argument as in the proof for 
$(ii) \Rightarrow (i)$, 
we can show that $P(N, \varphi_0)$ is valid for all $N \in {\bf N}$ and $\varphi _0 \in \SCR{S}({\bf R}) \backslash \{ 0 \}$ by induction. 
If we take $t=t_0$, we have $(ii)$ of Theorem \ref{T2}. \qed

\end{document}